\documentclass[11pt,english]{amsart}
\usepackage[letterpaper,top=2cm,bottom=2.5cm,left=3cm,right=3cm,marginparwidth=1.75cm]{geometry}
\usepackage{color}
\usepackage[colorlinks=true, allcolors=blue]{hyperref}

\usepackage{amssymb}
\usepackage{lineno}

\usepackage{cancel}

\usepackage{enumitem}

\usepackage{graphicx}
\usepackage{xcolor}
\usepackage{amsmath}
\usepackage[normalem]{ulem}

\usepackage{subcaption} % Necesario para el entorno subfigure

\newtheorem{theorem}{Theorem}[section]
\newtheorem{thmy}{Theorem}
\newtheorem{lemma}[theorem]{Lemma}
\newtheorem{proposition}[theorem]{Proposition}

\theoremstyle{definition}

\theoremstyle{remark}
\newtheorem{remark}[theorem]{Remark}

\numberwithin{equation}{section}

\usepackage{mathrsfs}
\usepackage{comment}

\newcommand{\bbZ}{\mathbb{Z}}
\newcommand{\bbC}{\mathbb{C}}

\newcommand{\bbN}{\mathbb{N}}

\newcommand{\Cc}{\widehat{\bbC}}

\newcommand{\fat}{\mathcal{F}}
\newcommand{\jul}{\mathcal{J}}

\newcommand{\cQ}{\mathcal{Q}}

\newcommand{\lau}{\lambda_1}
\newcommand{\lad}{\lambda_2}
\newcommand{\lat}{\lambda_3}

\newcommand{\LA}{\Lambda}
\newcommand{\la}{\lambda}
\newcommand{\TLA}{\mathbb{T}_\Lambda}
\newcommand{\WLA}{\wp_\Lambda}
\newcommand{\ALA}{\mathcal{A}(\Lambda)}
\newcommand{\fLA}{f_\Lambda}

\newcommand{\NLA}{N_{f,\Lambda}}

\DeclareMathOperator{\sing}{sing}
\DeclareMathOperator{\Int}{int}
\DeclareMathOperator{\Crit}{Crit}
\DeclareMathOperator{\Zeros}{Zeros}
\DeclareMathOperator{\Poles}{Poles}
\DeclareMathOperator{\Ext}{ext}
\DeclareMathOperator{\Fix}{Fix}

\begin{document}

\title{Dynamics of Newton's method for odd and even elliptic functions}

\author[Esparza-Amador]{Adri\'an Esparza-Amador}
\address{Universidad Austral de Chile}
\email{adrian.esparza@uach.cl}
%\thanks{The author was supported in part }

\author[Moreno Rocha]{M\'onica Moreno Rocha}
\address{Centro de Investigaci\'on en Matem\'aticas -- CIMAT, M\'exico.}
\email{mmoreno@cimat.mx}
%\thanks{The author was supported in part }

%\subjclass[2000]{Primary 37F10; Secondary 30D05, 37F50, 65H05}
%\date{\today}
%\dedicatory{}

\begin{abstract}
We investigate Newton's method applied to any odd or any even elliptic function with an arbitrary period lattice. For any function of this type whose set of poles coincides with its period lattice, we show that the Julia set of its Newton map is connected, as long as no Herman rings exist. Moreover, we provide sufficient conditions on the Newton's method of any odd or even elliptic function to exhibit wandering domains coexisting with attracting basins. This phenomenon was first reported by Florido and Fagella for the Newton's method applied to a family of entire functions; however, our approach does not employ the lifting technique. We also provide a detailed study of a one-parameter family of elliptic functions given by $\WLA+b$ with $\LA$ any triangular period lattice and $b\in\bbC$. We show that their associated Newton maps do not exhibit Herman rings or Baker domains, and any other Fatou component, including wandering, is bounded.

\end{abstract}

%\classification{37F10; 30D05, \textcolor{blue}{30D30, 33E05},  37F50}

\maketitle

\section{Introduction}\label{sec:Intro}
The study of the dynamical properties of root-finding methods applied to transcendental functions has attracted much attention in the holomorphic dynamics community in the past thirty-five years; see, for example, \cite{BFJK, BFJK2, Ber, Be93-2, BR}. %, Ho, Kn}.
Among some of the most important questions regarding these algorithms are the connectivity of the associated Julia set and the possible realisation of Herman rings, Baker and wandering domains which, by definition, contain points in the complex plane where the convergence of the algorithm fails.

For the case of the Newton-Raphson root-finding algorithm, which from now on we will refer to as the \emph{Newton's method}, the connectivity question has been answered in great generality for entire functions: Shishikura has shown that the Julia set of the Newton's method for a polynomial is always connected, \cite{S}. The same result has been extended by Bara\'nski,  Fagella, Jarque and Karpi\'nska for the case of transcendental entire functions; see \cite{BFJK} and its unified exposition in \cite{BFJK2}.
Regarding the existence of certain Fatou components where the convergence of the Newton's method fails, we highlight the recent work by Florido and Fagella in \cite{FF25}, which is close to our results. 
They study a family of transcendental periodic entire functions with zeros for which their Newton maps become meromorphic functions that can be projected to the cylinder under the exponential map. Then, using the lifting technique, they provide the first examples of Newton maps that exhibit Baker domains or wandering domains coexisting with attracting basins for the roots of the entire function. Naturally, the realisation of a Herman ring by a Newton map implies a disconnected Julia set. It remains an open question whether any root-finding algorithm can realise this type of Fatou component; see \cite{CGN}.

In this work, we study the dynamical properties of the iterative Newton's method applied to an elliptic function $\fLA:\bbC\to \Cc$ whose lattice of periods, $\LA\subset\bbC$, is arbitrary and fixed. 
Our interest resides primarily on the connectivity of the Julia set associated with the transcendental meromorphic function
\[N_{\fLA}(z)= z - \frac{\fLA(z)}{\fLA'(z)},\]
and the existence of Herman rings, Baker, or wandering domains.

Before proceeding to the description of our results, we point out that a study of the \emph{continuous version} of Newton's (damped) method for elliptic functions has been performed by Helminck, Twilt and their coauthors, providing sufficient conditions for the phase portraits of (desingularised) elliptic vector fields to remain stable; see \cite{HT} and references within. The relationship of their work with ours is confined to the local phase portraits at its zeros, poles, and critical points of the associated elliptic Newton flow, which directly translates into the standard concepts of stability for fixed points and multiplicity of poles for the corresponding discrete Newton's method. 

\subsection*{Statements of results}
Let $\LA\subset\bbC$ denote an arbitrary lattice and consider any elliptic function $\fLA$. Then the Newton map  $N_{\fLA}$ is a \emph{$\LA$-equivariant function}: that is, $N_{\fLA}$ commutes with the action of $\LA$ by translations. This property is at the core of most of our main results as we shall see next.

In what follows, we only consider the collection of all odd and all even elliptic functions with the same period lattice $\LA$, which we denote by $\ALA$. If $\fLA$ belongs to this set, then its Newton map is an odd function. This, together with the $\LA$-equivariance, implies that most Fatou components of $N_{\fLA}$ are bounded (refer
to  Proposition~\ref{prop:bounded} for a precise statement), while any half-period lattice point is a \emph{point of symmetry} for $N_{\fLA}$, see Proposition~\ref{prop:symm}. In particular, if $\fLA\in\ALA$ and its set of poles coincides with its period lattice (Lemma~\ref{lem:tipo_pol} provides a characterisation of these elliptic functions), then its Newton map $N_{\fLA}$ has a repelling fixed point at every  point $\la\in \LA$.

Our first main result deals with the connectivity of the Julia set, denoted by $\jul(N_{\fLA})$, for any elliptic function $\fLA\in\ALA$ whose set of poles is exactly its period lattice $\LA$. To our knowledge, this is the first connectivity result of its kind concerning the Newton's method of an elliptic function.

\begin{thmy}
    \label{thm:conn_gen}
    Fix an arbitrary lattice $\LA$ and let $N_{\fLA}$ denote the Newton's method of an odd or even elliptic function $\fLA$ for which $\Poles(\fLA)=\LA$. If $N_{\fLA}$ has no Herman rings, then its Julia set is connected. 
\end{thmy}
We briefly describe the strategy of the proof of Theorem~\ref{thm:conn_gen}: provided that the Fatou set, denoted by $\fat(N_{\fLA})$, contains no Herman rings, we show that any other Fatou component is simply connected. For otherwise, if $U\subset\fat(N_{\fLA})$ is multiply connected, the application of several results from \cite{BFJK2}  implies the existence of a \emph{weakly repelling} fixed point in a bounded complementary component of $U$. The hypothesis on the poles of $\fLA$ then shows that $U$ is unbounded and the last step reduces to proving that $U$ is a simply connected Baker domain or a simply connected escaping wandering domain, which is the essence of our second main result.

\begin{thmy}
\label{thm:unb_comp}
    Fix an arbitrary lattice $\LA$ and let $N_{f_\LA}$ denote the Newton's method of an odd or even elliptic function $f_\LA$ for which $\Poles(\fLA)=\LA$. 
    If $U$ %\subset\fat(N_{\fLA})$ is 
    is an unbounded Fatou component, then $U$ is a periodic Baker domain or an escaping wandering domain. Furthermore, $U$ is simply connected.
\end{thmy}

We now turn to discuss the realisation of Herman rings, Baker domains and wandering domains. The hypothesis in Theorem \ref{thm:conn_gen} on the absence of Herman rings is addressed for the Newton's method of a one-parameter family of (even) elliptic functions given by
\[\wp_{\LA, b}(z)=\WLA(z)+b,\]
where $b\in\bbC$ and $\LA$ is triangular, however other results for arbitrary lattices are also given. The Newton maps associated with $\wp_{\LA,b}$ are simply denoted by $N_b$. When $N_b$ is restricted to a fundamental parallelogram $\cQ$, it has four free, simple critical points unless $\LA$ is triangular; in this case, $N_b|\cQ$ has exactly two free, double critical points. 
 For any $q\geq 1$, we show that if a $q$-cycle of Herman rings is realised by $N_b$ and $\LA$ is arbitrary, then the cycle presents a \emph{nested configuration} and each Herman ring surrounds the same pole; see Proposition~\ref{prop:onepole}. And although Newton's method of an elliptic function $\fLA\in\ALA$ may realise unbounded components as established in Theorem~\ref{thm:unb_comp}, exploiting the properties of triangular lattices, we are able to prove our third main result.

\begin{thmy}\label{thm:NoHR}
    Let $\LA$ be a triangular lattice. Then, for any $b\in \bbC$, 
    the Newton map
    \[N_b(z) = z - \frac{\WLA(z)+b}{\WLA'(z)}\]
    does not have cycles of Herman rings of any period and every Fatou component of $N_b$ is bounded. Furthermore, its Julia set is connected.
\end{thmy}

Our final result is related to the examples constructed by Florido and Fagella in \cite{FF25} on the existence of wandering domains and (super)attracting basins of roots in the same dynamical plane of $N_{\fLA}$ and $\fLA\in\ALA$ with $\LA$ an arbitrary lattice.

\begin{thmy}\label{thm:wandering}
For any given lattice $\LA$, let $N_{\fLA}$ denote the Newton's method for an odd or even elliptic function $\fLA$. Consider the function $F_{\la}(z)=N_{\fLA}(z)+\la$ for some fixed $\la\in \LA^*$ and let $q\geq 1$. If $F_{\la}$ has a $q$-cycle of (super)attracting  basins containing only free critical points,  %in $\Zeros(\fLA'')\setminus \Zeros(\fLA)$},
then $N_{\fLA}$ has (at least) $2q$ distinct orbits of escaping wandering domains that coexist with the invariant (super)attracting basins for the roots of $\fLA$.
%For any given lattice $\LA$, let $N_{\fLA}$ denote the Newton's method for an odd or an even elliptic function $\fLA$. Consider the function $F_{\la}(z)=N_{\fLA}(z)+\la$ for some fixed $\la\in \LA^*$. If $F_{\la}$ has a (super)attracting $q$-cycle for some $q\geq 1$, then $N_{\fLA}$ has (at least) $2q$ distinct orbits of simply connected wandering domain. 
\end{thmy}

In contrast with the lifting method used in \cite{FF25} to generate wandering domains, the proof of Theorem~\ref{thm:wandering} rests on the $\LA$-invariance property of $N_{\fLA}$ and Tsantaris result in \cite[Theorem 1.5]{Tsa}, which shows that any two commuting transcendental meromorphic functions with infinitely many prepoles share the same Julia set. Surprisingly, the dynamical plane of Newton maps $N_b$ may also exhibit (bounded) wandering domains together with the (super)attracting basins of the roots of $\WLA+b$; see Proposition~\ref{ex:wanderN_b} where we describe parameters that give rise to pairs of distinct orbits of wandering domains. Theorem~\ref{thm:NoHR} implies these wandering domains are simply connected.

The proofs of Theorem~\ref{thm:conn_gen} and Theorem~\ref{thm:unb_comp} are given in Section~\ref{sec::Conn}, while in Section~\ref{sec:ParmFam}, the proof of Theorem~\ref{thm:NoHR} is presented. 
 Section~\ref{sec:wan} contains the proofs of Theorem~\ref{thm:wandering} and Proposition~\ref{ex:wanderN_b}, including some numerical experiments.

\section{Preliminaries}\label{sec:prelim}
We provide a brief overview of the theory of elliptic functions and then present the relevant features of the Newton's method of an elliptic function. We refer the reader to the classical works of Du Val \cite{DuVal} for further details of elliptic functions and Bergweiler \cite{Ber} for the dynamics of transcendental meromophic functions. We also recommend the modern presentation by Kotus \& Urba\'nski in \cite{Kotus} on the theory of elliptic functions and their dynamics.

\subsection{Elliptic functions}\label{subsec:elliptic}

Throughout this work we consider a fixed lattice
\[\LA=\{\la:\la=m\lau+n\lad, m,n\in \bbZ\}\]
whose generators, $\lau, \lad\in \bbC$, satisfy $\text{Im}(\lad/\lau)>0$. Each lattice $\LA$ determines a pair of \textit{invariants},
\[g_2:=g_2(\LA)=60\sum_{\lambda\in\LA\setminus\{0\}}\lambda^{-4}\qquad\text{and}\qquad g_3:=g_3(\LA)=140\sum_{\lambda\in\LA\setminus\{0\}}\lambda^{-6}\]
that satisfy $g_2^3-27g_3^2\neq 0$, and vice versa (see Corollary 16.5.11 in \cite{Kotus}). For further reference, we say that $\LA$ is \textit{triangular} if and only if $g_2=0$ and $g_3\neq 0$, \cite{DuVal}.

$\LA$ acts by translations in $\bbC$, so we denote by $\cQ_0$ the fundamental parallelogram with vertices in $0,\lau, \lad$ and $\lat=\lau+\lad$. We will refer to $\cQ_0$ as the \textit{principal} parallelogram.

The {\it Weierstrass elliptic function}, $\WLA$, and its derivative, $\WLA'$, are related by the differential equation 
\begin{equation}
    \label{eq:rel_p-pprime}
    \WLA'(z)^2 = 4\WLA(z)^3 - g_2\WLA(z) - g_3.
\end{equation}
From this equation, one can see that $\WLA$ has three distinct residue classes of critical points determined by the \textit{half-periods} $\omega_i=\la_i/2$ for $i=1,2,3$. The \textit{critical values} of $\WLA$, customarily denoted as $e_i:=e_i(\LA)=\WLA(\omega_i)$, are all distinct. For further reference, we will denote the set of half-periods of $\WLA$ by
\[\frac 12 \LA \setminus \LA := \{\omega_1, \omega_2, \omega_3\}+\LA,\]
and remark that $\frac 12 \LA$ naturally contains $\LA$.

If $\fLA$ denotes any other elliptic function with period lattice $\LA$, then it belongs to the rational field $\bbC(\WLA, \WLA')$.

\begin{theorem}
    \label{thm:all_ellip}
    Every elliptic function $\fLA:\bbC\to \Cc$ with period lattice $\LA$ can be written as $\fLA(z) = R(\WLA(z))+\WLA'(z)Q(\WLA(z))$, where $R$ and $Q$ are rational functions with complex coefficients. The converse is also true; that is, every function $\fLA$ of this form is elliptic. 
\end{theorem}

Denote by $o_f\geq 2$ the order of $f=\fLA$. An elliptic function does not have omitted values (as $\fLA$ takes each value $z\in \bbC$ exactly $o_f$ times) and, due to its double periodicity, it does not have asymptotic values. 
For later use, we characterise those elliptic functions with poles given exactly by the period lattice $\LA$.

\begin{lemma}
    \label{lem:tipo_pol}
    Let $\fLA$ be an arbitrary elliptic function with period lattice $\LA$. Then, $\Poles(\fLA)=\LA$ if and only if $\fLA(z)=P(\WLA)+\WLA'\cdot S(\WLA)$, for some polynomials $P$ and $S$ with complex coefficients. 
\end{lemma}

\subsection{The Newton method's of an elliptic function}\label{subsec::Newton}

For a fixed lattice $\LA$, we study the dynamics  of the Newton's method of an (uni-parametric) elliptic function $f=\fLA$, denoted by
\begin{equation}\label{eq:NE}
    N_{f,\LA}(z)=z-\frac{\fLA(z)}{\fLA'(z)}.
%N_f(z)=z-\frac{f(z)}{f'(z)}.
\end{equation}
For simplicity, the dependence of the elliptic function and its Newton's method on the lattice $\LA$ will be suppressed most of the times.

The transcendental meromorphic function $N_f:\bbC\to \Cc$ has an essential singularity at infinity, with no omitted values or (finite) asymptotic values. This function is no longer double periodic, but commutes with the action of $\LA$ by translations. Since this property will be extensively used throughout this work, it will be referred to as the $\LA$-{\it equivariance property}: that is,
\begin{equation}\label{eq:conmuta}
    N_f(z+\la) =N_f(z)+\la,\qquad \text{for all}~\la\in \LA, \text{ and all } z\in \bbC.
\end{equation}

An immediate consequence of the $\LA$-equivariance is the following.

\begin{lemma}\label{lem::inv-LA}
    Fix a lattice $\LA$ and consider the Newton map $N_{f}(z)= z-\frac{f(z)}{f'(z)}$, where $f$ is a (nonconstant) elliptic function. Then, $\jul(N_{f})=\jul(N_{f})+\LA$ and in consequenece, $\fat(N_{f})=\fat(N_{f})+\LA$. 
\end{lemma}

As for the case of Newton's method for rational functions, a zero of the elliptic function $f$ is an attracting fixed point of $N_f$ and in particular, simple zeros of $f$ become super-attracting. The following result is easily obtained from (\ref{eq:NE}).

\begin{lemma}\label{lem:propN}
Let $f:\bbC\to \Cc$ be an elliptic function with respect to $\LA$ of order $o_{f}\geq 2$ and let $N_f$ denote its Newton's method. Then, 
\begin{enumerate}
    \item a zero of $f$ of order $m\geq 1$ is a (super)attracting fixed point of $N_f$ with multiplier $\la=\frac{m-1}{m}$,
    while a pole of $f$ of order $m\geq 1$ is a repelling fixed point of $N_f$ with multiplier $\la=\frac{m+1}{m}$. In particular,
\[\Fix(N_f)=\Zeros(f) \sqcup \Poles(f).\]

    \item The set of poles of $N_f$ consists of the zeros of $f'$ that are not zeros of $f$, that is,
\[\Poles(N_f)=\Zeros(f')\setminus \Zeros(f) .\]
    \item The critical points of $N_f$ are given by the zeros of
\begin{equation}\label{eq:NPrima}
    N_f'(z) = \frac{f(z)f''(z)}{(f'(z))^2}
\end{equation}
and all their $\LA$-orbits. 
\end{enumerate}
\end{lemma}

Following \cite{NP}, we say that a critical point of $N_f$ is {\it additional} if it is a zero of $f''$, and it is {\it free} if the critical point is not fixed under $N_f$. It follows that each free critical point is additional, too.

\begin{remark}\label{rem:equivariance}
The $\LA$-equivariance of $N_f$ has two more immediate consequences. First, by (3) in the previous lemma, its set of critical values accumulate at $\infty$, which is itself an asymptotic value of $N_f$, so in principle, $N_f$ may realise Baker domains or wandering domains.

Second, the equivariance property allows us to quotient the action of $N_f$ into the flat $2$-torus $\TLA$, inducing a new function $n_f:\TLA^*\to\TLA$, that satisfies
\begin{equation}\label{def:proj_map}
n_f\circ \Pi(z)=\Pi\circ N_f(z),\qquad \forall z\in \bbC\setminus\Poles(N_f),
\end{equation}
where $\Pi:\bbC\to \TLA$ is the natural projection and $\TLA^*=\TLA\setminus\Pi(\Poles(N_f))$. Observe that $n_f$ is a meromorphic function of finite type over $\TLA^*$ (see also Section 16.2 in \cite{Kotus}).
\end{remark}

The Newton's method of an odd or an even elliptic function has further elementary properties which are described in the next two results.

\begin{proposition}
    \label{prop:N_odd}
    If $f$ is an odd or even elliptic function, then $N_f$ is an odd function. 
\end{proposition}
\begin{proof}
    ~It follows from Theorem \ref{thm:all_ellip} that an elliptic function $f$ is even (respectively odd) if and only if the rational function $Q$ is identically zero (respectively the rational function $R$ is constant). Hence, the derivative map interchanges the parity of the function $f$ and the quotient $f/f'$ is always odd. In this way, the Newton map is the sum of odd functions. 
\end{proof}

The Newton's method of an odd or even elliptic function satisfies a more general symmetry among their orbits:  Given three distinct points $p,z,w\in \bbC$, we say that $z$ and $w$ are {\it symmetric with respect to $p$} if
\[p-z= w-p.\]

\begin{proposition}[Symmetric orbits]
\label{prop:symm}
Let $f$ be an odd or even elliptic function with period lattice $\LA$ and denote by $N_f$ its Newton's method. If $z,w\in \bbC$ are two distinct points that are symmetric with respect to a half-period $p\in\frac 12\LA\setminus\LA$, then the positive orbits of $z$ and $w$ are also symmetric with respect to $p$. In other words, for any integer $k\geq 0$, the relation
\[p-N_f^k(z)=N_f^k(w)-p\]
holds, and in consequence $N_f^k(w)=2p-N_f^k(z)$ for all $k\geq 0$.

%In particular, if $\la=0$ it follows that $N_f^k(-z)=-N_f^k(z)$.
\end{proposition}
\begin{proof}
~The relation $p-z=w-p$ is equivalent to $w=2p-z$ with $2p=\la\in\LA$. From Proposition~\ref{prop:N_odd} and the $\LA$-equivariance property in (\ref{eq:conmuta}) one obtains
\[N_f(w)=N_f(\la-z)=-N_f(z-\la)=-(N_f(z)-\la)=\la-N_f(z),\]
and the desire relation follows directly.

\end{proof}

\begin{lemma}\label{ex:NewtonP}%\footnote{esto deber\'ia ser un lema o un  -- voy más por un lemma, ya que se menciona en una de las demos}
     Fix an arbitrary lattice $\LA$ and consider the Newton's method for the Weierstrass $\WLA$ function.
     Then, the additional critical points of $N_{\WLA}$ are completely determined by the invariant $g_2=g_2(\LA)$. In particular, if $\LA$ is 
     triangular, then $N_{\WLA}$ has no additional critical points and the Julia set $\jul(N_{\WLA})$ is connected.
\end{lemma}
\begin{proof}
    ~From Lemma \ref{lem:propN} (3), we know that 
    \[
    N_{\WLA} '(z) = \frac{\WLA(z)\cdot\WLA''(z)}{(\WLA'(z))^2}. 
    \]
    By deriving both sides of (\ref{eq:rel_p-pprime}), we obtain a similar relation between the Weierstrass elliptic function and its second derivative 
    \begin{equation}
        \label{eq:second_der}
        \WLA''(z) = 6\WLA(z)^2-\frac{1}{2}g_2.
    \end{equation}
    From this relation, we have a full description of the critical points, namely 
    \begin{equation}
        \Crit(N_{\WLA}) = \Zeros(\WLA)\cup \left\{\WLA^{-1}\left(\pm\sqrt{\dfrac{g_2}{12}}\right)\right\}.
    \end{equation}
    
 Therefore, additional critical points are completely determined in terms of the lattice invariant $g_2$.  
    In the particular case where $g_2=g_2(\LA)=0$ (and consequently, $\LA$ is triangular) 
    the critical set of $N_{\WLA}$ is reduced to the roots of $\WLA$ 
 which become super-attracting fixed points of $N_{\WLA}$ 
 and thus, the Fatou set $\fat(N_{\WLA})$ consists only of 
 super-attracting basins. 
 Furthermore, a suitable application of the Riemann-Hurwitz formula indicates that each Fatou component is simply connected and therefore, the Julia set $\jul(N_{\WLA})$ is connected.
\end{proof}

\section{Connectivity of the Julia set of elliptic Newton's method}\label{sec::Conn}
This section contains the proof of our first main result.

\vspace{0.3cm}
\noindent{\bf Theorem A.}     Fix an arbitrary lattice $\LA$ and let $N_{\fLA}$ denote the Newton's method of an odd or even elliptic function $\fLA$ for which $\Poles(\fLA)=\LA$. If $N_{\fLA}$ has no Herman rings, then its Julia set is connected.\\

The strategy of the proof is standard as we show that, provided that $\fat(N_f)$ contains no Herman rings, then any other Fatou component is simply connected, for if, otherwise, a Fatou component is multiple connected (and not a Herman ring), then we employ two lemmas from Baranski et al. in \cite{BFJK2} to obtain the existence of weakly repelling fixed points that, in turn, will yield a contradiction. 

Recall that a function $g$ has a weakly repelling fixed point at $z_0$ if $z_0=g(z_0)$ and $|g'(z_0)|>1$ or $g'(z_0)=1$.
For completeness, we provide the statements of the aforementioned lemmas below and introduce the required notation. If $X\subset \bbC$ is a compact set, we denote by $\Ext (X)$ the connected component of $\widehat{\bbC}\setminus X$ that contains the infinity. We set $K(X)=\widehat{\bbC}\setminus \Ext(X)$, which is a closed and bounded set. From the maximum modulus principle, it follows that if $f$ is a holomorphic map in a neighbourhood of $K(X)$, then $f(K(X))=K(f(X))$. Finally, if $\gamma\subset \bbC$ is a Jordan curve, $\Int(\gamma)$ will denote the bounded component of $\bbC\setminus \gamma$. 

\begin{lemma}[Lemma 2.4 in \cite{BFJK2}]
    \label{lem:polesloop}
    Let $f:\bbC\to\widehat{\bbC}$ be a meromorphic transcendental map or a rational map for which infinity belongs to the Julia set. Let $\gamma\subset\bbC$ be a closed curve in a Fatou component $U$ of $f$, such that $K(\gamma)\cap\jul(f)\neq\emptyset$. Then there exists $n\geq0$, such that $K(f^n(\gamma))$ contains a pole of $f$. Consequently, if $U$ is multiply connected then there exists a bounded component of $\Cc\setminus f^n(U)$, which contains a pole. 
\end{lemma}

\begin{lemma}[Lemma 2.8 in \cite{BFJK2}]
    \label{lem:boudariesout}
    Let $\Omega\subset\bbC$ be a bounded domain with finite Euler characteristic and let $f$ be a meromorphic map in a neighborhood of $\overline{\Omega}$. Assume that there exists a component $D$ of $\widehat{\bbC}\setminus f(\partial \Omega)$, such that: 
    \begin{itemize}
        \item[(a)] $\overline{\Omega}\subset D$, 
        \item[(b)] there exists $z_0\in\Omega$ such that $f(z_0)\in D$.  
    \end{itemize}
    Then $f$ has a weakly repelling fixed point in $\Omega$. 
\end{lemma}

We start the proof of Theorem \ref{thm:conn_gen} by determining in the next proposition the Fatou components of $N_f$ that are bounded. Then, we divide the proof in two cases: Section~\ref{sec:Bd} addresses the connectivity of bounded Fatou domains, while unbounded Fatou components (either preperiodic or wandering) are discussed in Section~\ref{sec:Ubd}.

\begin{proposition}\label{prop:bounded}
Fix  a lattice $\LA$ and let $N_f$ be the Newton's method of an odd or even elliptic function $f$. If $U\subset \fat(N_f)$ is a component that eventually maps into, or belongs to, a cycle of (super)attracting, parabolic, Siegel disks or Herman rings, then $U$ lies within a fundamental domain for the action of $\LA$ in $\bbC$ and, therefore, $U$ is bounded.

Then, if $N_f$ has no Baker domains or wandering domains, each component of $\fat(N_f)$ lies within a fundamental domain for the action of $\LA$ in $\bbC$ and, therefore, it is bounded. 
\end{proposition}%}
\begin{proof}
    ~Let $U$ be a Fatou component and assume that there exist points $z,w\in U$ and $\la\in \LA$ such that $w=z+\la$. We will show that $\la=0$.

    Assume first that $U$ is a $q$-periodic Fatou component for some $q\geq 1$.
    \begin{enumerate}
        \item Consider the case where the limit functions of $\{N_f^k\}_{k\geq 1}$ are all constant throughout the $q$-cycle: then, there exists $z_0\in \overline{U}\cap \bbC$ such that $N_f^{kq}\to z_0$ locally uniformly throughout $U$. 
        At the same time, the $\LA$-equivariance property in (\ref{eq:conmuta}) implies $N_f^{kq}(w)=N_f^{kq}(z)+\la\to z_0+\la$, and the uniqueness of the limit implies $\la=0$.   
        \item Consider the case where $U$ is a Siegel disc and neither $z$ nor $w$ are the $q$-periodic Siegel point $z_0\in U$. Then there exist simple closed curves, $\gamma_z$ and $\gamma_w$ in the foliation of $U$, where the orbits of $z$ and $w$ accumulate under $\NLA^{kq}$, respectively. Using the $\LA$-equivariance property, the orbit of $w$ also accumulates in $\gamma_z+\la$, therefore $\gamma_w=\gamma_z+\la\subset U$. 
        Since $\gamma_z$ winds around $z_0$, then $\gamma_z+\la$ winds around $z_0+\la$ and the simple connectivity of $U$ implies that $z_0+\la$ also lies in $U$. Then again, $\LA$-equivariance shows that $z_0+\la$ is fixed under $N_f^q$, which can only happen when $\la=0$. The case where either $z$ or $w$ is the $q$-periodic Siegel point in $U$ follows immediately.
        \item Assume $U$ is a Herman ring. Since every iterate $N_f^j(U)$ is a Herman ring containing points $N_f^j(z)$ and $N_f^j(z)+\la$, we can assume without loss of generality that the bounded component of $\bbC\setminus U$ contains a pole $p\in\Poles(N_f)$. Arguing as in the case of Siegel discs, there are simple closed curves $\gamma_z$ and $\gamma_w=\gamma_z+\la$ in the foliation of $U$ winding around both $p$ and $p+\la$. Then again, $\gamma_z+\la$ (and hence $\gamma_z$) must also wind around $p+\la$ and $p+2\la$. Continuing in this way, one sees that $\gamma_z$ winds around $p+k\la$, for $k\geq 1$.  Although $\partial U$ may contain the point at infinity, each foliation curve is bounded, implying that $\la=0$.
    \end{enumerate}
  To finish the proof, assume $U$ is a preperiodic Fatou component that contains $z$ and $w$ as before. Then there exists an iterate $k\geq 1$ so that $N_f^k(z)$ and $N_f^k(w)=N_f^k(z)+\la$ belong to a periodic Fatou component. The previous cases imply that $\la=0$.
\end{proof}

\subsection{Connectivity of bounded domains}\label{sec:Bd} 
We begin with a technical result that establishes the winding number of the iterates under $N_f$ of a simple closed curve symmetric with respect to a point in $\frac 12 \LA$. Observe that no restriction is imposed on the set $\Poles(f)$. 

\begin{lemma}
    \label{lem:simm_curves} Fix a lattice $\LA$ and let $N_f$ be the Newton's method of an even or odd elliptic function $f$.
Let $z_0\in\frac{1}{2}\LA$ be a period or half-period point of $f$. 
    Suppose $\gamma\subset\bbC$ is a simple closed curve such that $z_0\in\Int(\gamma)$. If $\gamma$ is symmetric with respect to $z_0$, then $z_0\in \Int(N_f^n(\gamma))$ for $n\geq1$. %\footnote{generalizado}
\end{lemma}
\begin{proof}
    ~By induction, since $N_f^{n+1} = N_f\circ N_f^n$, it is enough to prove the case for $n=1$. 

    First, suppose $z_0=\la'\in\LA$ is a lattice point. From the $\LA$-equivariance property (\ref{eq:conmuta}), we may assume $\lambda'=0$. Hence, the curve $\gamma$ is symmetric with respect to the origin. Since $f$ is odd or even, it follows from Proposition \ref{prop:N_odd} that $N_f$ is an odd function, which implies that 
    \[ N_f^n(-z) = -N_f^n(z),\]
    for all $n\geq0$, and for all $z\in\bbC$, whenever the iterates are well-defined. In other words, the symmetry with respect to the origin is preserved under the map $N_f$. Hence, $N_f(\gamma)$ is symmetric with respect the origin and thus  $0\in\Int(N_f(\gamma))$. 

    Now, if $z_0=\omega'\in\frac 12\LA\setminus\LA$ is a half-period point. Analogous to the previous case, from Proposition \ref{prop:symm} the symmetry with respect to the half-period $\omega'$ is preserved under $N_f$, which implies $\omega '\in\Int(N_f(\gamma))$, and the proof is done.
    
\end{proof}

The next result establishes the connectivity of the Julia set $\jul(N_f)$ under the hypotheses of Theorem \ref{thm:conn_gen} and in the case when each connected component of the Fatou set is bounded (in other words, there are no Baker domains or unbounded wandering domains).

\begin{lemma}
    \label{lem:simplybounded}
    Fix a lattice $\LA$ and let $N_f$ be the Newton's method of an odd or even elliptic function $f$ for which $\Poles(f)=\LA$. If $N_f$ has no Herman rings or unbounded Fatou components, then each component of $\fat(N_f)$ is simply connected. Equivalently, the Julia set $\jul(N_f)$ is connected. 
\end{lemma}

\begin{proof}
    ~First of all, note that, from Proposition \ref{prop:bounded}, every Fatou component is contained within a fundamental domain. 
    
    Let $U\subset\fat(N_f)$ be a multiply connected component. Since $N_f$ has no Herman rings, it follows that $U$ is a (pre)periodic component of a (super)attracting or parabolic basin, or a pre-periodic component of a cycle of Siegel discs, or a wandering domain (see Theorem \ref{thm:wandering}). 

    From Lemma \ref{lem:polesloop} there exists $n_0\geq0$ such that a bounded component of $\Cc\setminus(N_f^{n_0}(U))$ contains a pole of $N_f$, that is, there is a simple closed curve $\gamma\subset N_f^{n_0}(U)$ such that $p\in\Int(\gamma)$ for some pole $p$ of $N_f$. Since $N_f^{n_0}(U)\subset\fat(N_f)$ by invariance,  we may suppose, without loss of generality, that $n_0=0$. Then $\gamma\subset U$ with $p\in\Int(\gamma)$. 

    Recall that, in the case that $U$ is a (pre)periodic component of a (super)attracting or parabolic basin, or a pre-periodic component of a cycle of Siegel discs, then every curve in $U$ is eventually contractible. It is not difficult to see that we may consider our curve $\gamma\subset U$ such that $\Int(\gamma)\cap\Int(N_f(\gamma))=\emptyset$. This last statement is also valid when $U$ is a wandering domain. 

    In this way, taking $\Omega=\Int(\gamma)$, $D=\Ext(N_f(\gamma))$, and $z_0=p$, Lemma \ref{lem:boudariesout} implies the existence of a weakly repelling fixed point $w_0$ contained in $\Omega$. Since $\Poles(f)=\LA$, from Lemma \ref{lem:propN} we have that the weakly repelling fixed point $w_0=\la'$ must be a lattice point. 

    We claim that $\gamma\subset U$ can be chosen to be symmetric with respect to $\la'$: If $\gamma$ were not symmetric, consider $\gamma'$ as the reflection of $\gamma$ with respect to the point $\la'$. Then, by construction, the curve $\gamma''=\partial K(\gamma\cup\gamma')$ is a simple closed curve, symmetric with respect to $\la'$. 

    Finally, we have a simple closed curve, symmetric with respect to a lattice point $\la'$, such that $\la'\in\Int(\gamma)$. From Lemma \ref{lem:simm_curves}, $\la'\in\Int(N_f^n(\gamma))$ for $n\geq1$, which contradicts the election of $\gamma$.
    We conclude that $U$ must be simply connected. 
\end{proof}

\subsection{Connectivity of unbounded components}\label{sec:Ubd}
As stated in the preliminaries, the Newton's method of any given elliptic function has an asymptotic value at infinity, so in principle it may realise (unbounded) wandering domains or Baker domains. For elliptic functions in the set $\ALA$, the next result shows that these are the only possible unbounded Fatou components of its Newton's method and they must be simply connected.
\\

\noindent{\bf Theorem B.}
    Fix an arbitrary lattice $\LA$ and let $N_{f_\LA}$ denote the Newton's method of an odd or even elliptic function $f_\LA$ for which $\Poles(\fLA)=\LA$. 
    If $U$ %\subset\fat(N_{\fLA})$ is 
    is an unbounded Fatou component, then $U$ is a periodic Baker domain or an escaping wandering domain. Furthermore, $U$ is simply connected.

\begin{proof}
    ~If $U$ is an unbounded, open and connected component of the complex plane, then there exists $z\in U$ and $\la\in \LA^*$ so that $z, z+\la\in U$. Proposition~\ref{prop:bounded} implies that $U$ is either periodic or preperiodic to a $q$-cycle of Baker domains or $U$ is a wandering domain. In both cases, there exists $z_0\in\Cc$ and an increasing subsequence $k_j\to \infty$, so that $N_f^{k_j}\to z_0$ locally uniformly throughout $U$ as $k_j\to\infty$. The $\LA$-invariance property implies that, as $z+\la\in U$,
    \[
    N_f^{k_j}(z+\la) = N_f^{k_j}(z)+\la \to z_0+\la,\ k_j\to\infty,
    \]
    which is impossible unless $z_0=\infty$. Since $N_f$ has no finite asymptotic values, we conclude that either $U$ belongs to a $q$-cycle of Baker domains (associated to a unique Baker point at infinity) or $U$ is a wandering domain with infinity as its unique constant limit function, therefore it is escaping.

Now, suppose $U$ is a multiple connected domain. From Lemma \ref{lem:polesloop} we may assume there exists a simply closed curve $\gamma\subset U$ and some $p\in\Poles(N_f)$ so that $p\in\Int(\gamma)$. Proceeding as in the proof of Lemma \ref{lem:simplybounded}, we may select $\gamma$ such that
\[\Int(\gamma)\cap\Int(N_f(\gamma))=\emptyset.\]
Indeed, the case when $U$ is an unbounded wandering domain follows directly, while the case when $U$ is a Baker domain is a consequence of the existence of absorbing regions (see \cite[Theorem 2.1]{BFJK2}), implying that every closed curve in $U$ is eventually contractible under iteration. Then, Lemma \ref{lem:boudariesout} implies the existence of a repelling fixed point $\la'\in\LA$ for $N_f$ inside $\Int(\gamma)$. A symmetrisation of $\gamma$ with respect to $\la'$ yields a contradiction. 

\end{proof}

The proof of Theorem~\ref{thm:conn_gen} is now a direct consequence of Lemma~\ref{lem:simplybounded} and Theorem~\ref{thm:unb_comp}.

\section{Dynamics of a parametric family of Newton maps}%Connectivity for the triangular case}
\label{sec:ParmFam}

Throughout this section, we consider a one-parameter family of elliptic functions given by
\[\wp_{\LA, b}(z)=\WLA(z)+b,\qquad b\in \bbC,\]
and analyse the Newton's method of each element of the parametric family, namely
\[N_b(z):=N_{\WLA,b}(z)=z-\frac{\WLA(z)+b}{\WLA'(z)}.\]

We begin with a discussion of the critical points of $N_{b}$ over an arbitrary period lattice with emphasis on triangular period lattices, which constitute one of the main hypotheses in Theorem~\ref{thm:NoHR}, which is recalled below for convenience. Its proof will be a consequence of the auxiliary results.

\vspace{0.3cm}

\noindent{\bf Theorem C.}
    Let $\LA$ be a triangular lattice. Then, for any $b\in \bbC$,   the Newton map
    \[N_b(z) = z - \frac{\WLA(z)+b}{\WLA'(z)}\]
    does not have cycles of Herman rings of any period and every Fatou component of $N_b$ is bounded. Furthermore, its Julia set is connected.

\vspace{0.3cm}

Let $\LA$ be an arbitrary lattice with invariants $g_2,g_3\in \bbC, g_2^3-27g_3^2\neq 0$, and recall that $e_1=e_1(\LA), e_2=e_2(\LA)$ and $e_3=e_3(\LA)$ are the critical values of $\WLA$. 
The critical set for the Newton map $N_{b}$ becomes
\begin{equation}
    \label{eq:crit_}
    \Crit(N_{b})=\left\{\WLA^{-1}(-b)\right\}\cup \left\{\WLA^{-1}\left(\pm\sqrt{\frac{g_2}{12}} \right)\right\}.
\end{equation}
First, observe that since $\WLA$ has order 2, then the number of additional critical points (without counting multiplicity) of $N_{b}$ over any fundamental domain is either $2$ (if $g_2=0$) or $4$ (if $g_2\neq 0$). For critical points associated to the roots of $\WLA+b$, consider the case when
\[b\in \{-e_1, -e_2, -e_3\}.\]
The elliptic function $\WLA-e_i$ has a double zero in the half-periods $\omega_i+\LA$, 
so in particular, $N_{-e_i}$ has attracting fixed points at each point in this residual class. 

\begin{lemma}
\label{lemma:crittrian}
Let $\LA$ be an arbitrary lattice and let $b\in \bbC$ be such that $b\notin \{0,-e_1,-e_2,-e_3\}$. If $g_2(\LA)\neq0$ then each critical point of $N_{b}$ is simple, while if  $g_2(\LA)=0$, then each point in $\WLA^{-1}(0)$ is a double critical point of $N_{b}$.
\end{lemma}

\begin{proof}
    ~A direct computation using (\ref{eq:NPrima}) shows that
    \begin{equation}
        \label{eq:NBiprima}
        N_{b}''=\frac{(\WLA')^3\WLA''+(\WLA+b)(\WLA')^2\WLA'''-2(\WLA+b)\WLA'(\WLA'')^2}{(\WLA')^4},
    \end{equation}
    where we have omitted the argument on each function above for better readability. Also, one can easily obtain the expressions of the second, third and fourth derivatives of $\WLA$ over arbitrary lattices, given by
    \begin{equation*}
        \WLA''=6(\WLA)^2-g_2/2,\qquad \WLA'''=12\WLA \WLA',\qquad \WLA^{(IV)}=12(\WLA')^2+12\WLA \WLA''.
    \end{equation*}

    Assume $\zeta\in \WLA^{-1}(-b)$ for $-b\notin \{0,e_1,e_2,e_3\}$, so in particular $\WLA'(\zeta)\neq 0,\infty$. From (\ref{eq:NPrima}), $\zeta$ is a super-attracting fixed point of $N_{b}$. Moreover,
    \[N_{b}''(\zeta)=\frac{(\WLA'(\zeta))^3\WLA''(\zeta)}{(\WLA'(\zeta))^4}=\frac{\WLA''(\zeta)}{\WLA'(\zeta)}=\frac{6b^2}{\WLA'(\zeta)}\neq 0,\infty,\]
    that is, each critical point in $\WLA^{-1}(-b)$ is simple for the Newton map $N_{b}$. Suppose $\zeta\in\WLA^{-1}\left(\pm\sqrt{\frac{g_2}{12}} \right)$ for $g_2=g_2(\LA)\neq 0$. From (\ref{eq:NPrima}) and (\ref{eq:crit_}) we obtain $\WLA''(\zeta)=0$. Then equation (\ref{eq:NBiprima}) reduces to 
    \[N_{b}''(\zeta)=\frac{(\WLA(\zeta)+b)\WLA'''(\zeta)}{(\WLA'(\zeta))^2}=\frac{12(\WLA(\zeta)+b)\WLA(\zeta)}{(\WLA'(\zeta))}\neq 0,\infty,
    \]
    and $\zeta$ is a simple critical point, too.

    Now, if $g_2=0$ (so $\LA$ is triangular) take $\zeta\in \WLA^{-1}(0)$. On the one hand, $\WLA(\zeta)=0$ which implies $\WLA''(\zeta)=\WLA'''(\zeta)=0$. On the other hand, from the differential equation (\ref{eq:rel_p-pprime}) and the choice of a triangular lattice we obtain $(\WLA'(\zeta))^2=-g_3\neq 0$. From a direct inspection of (\ref{eq:NBiprima}), each term in the numerator of $N_{b}''(\zeta)$ becomes zero. The condition $(\WLA'(\zeta))^2=-g_3$ implies that $\WLA^{(IV)}(\zeta)=-12 g_3\neq 0$. Hence, taking the derivative of $N_{b}''$ and evaluating at $\zeta$, the unique non-zero term is
    \[N_{b}'''(\zeta)=\frac{(\WLA'(\zeta))^2(\WLA(\zeta)+b)\WLA^{(IV)}(\zeta)}{(\WLA'(\zeta)^4}=\frac{12b(\WLA'(\zeta))^4}{(\WLA'(\zeta))^4}=12 b\neq 0.\]
    We conclude that each critical point of $N_{b}$ in $\WLA^{-1}(0)$ is a double critical point.
\end{proof}

\begin{lemma}\label{lem:bees}
Let $\LA$ be an arbitrary lattice and consider the elliptic function $\WLA+b$ for $b=-e_i$, with $i=1,2,3$. If every free critical point of $N_{-e_i}$ lies in the immediate basins of $\{\omega_i\}+\LA$ then the Julia set of $N_{-e_i}$ is connected. In particular, if $\LA$ is triangular, then the Julia set of $N_{b}$ is connected for each $b\in\{0,-e_1, -e_2, -e_3\}$.
\end{lemma}
\begin{proof}
    ~We have already see that $\jul(N_0)$ is connected if $\LA$ is triangular. Let $\LA$ be an arbitrary period lattice.
    From the $\LA$-equivariance, it suffices to work with $b=-e_3$. Let $\mathcal{B}$ 
    denote the immediate basin of $\omega_3+\la$ (for some $\la\in\LA)$ which lies inside a period parallelogram $\mathcal{Q}$. Then $\mathcal{B}$ must contain a free critical point, say $c\in \WLA^{-1}(\pm\sqrt{g_2/12})\cap \mathcal{Q}$, and from Proposition~\ref{prop:symm}, it also contains its symmetric image $c'=2(\omega_3+\la)-c\in\mathcal Q$, which is also a free critical point.

    If $\LA$ is triangular (so $g_2=0$), the set of free critical points becomes
    \[\{\WLA^{-1}(0)\}=\{c,c'\}+\LA,\]
    where both are double critical points by Lemma~\ref{lemma:crittrian}.      
     If $g_2\neq 0$, denote by $d,d'$ the other two symmetrical and free critical points contained in $\mathcal{Q}$, and again from Lemma \ref{lemma:crittrian}, all free critical points are simple. If $d,d'\in \mathcal{B}$ then, regardless of the value of $g_2$, the Riemann-Hurwitz formula implies the equation
\begin{equation*}
(s-1)(m-2)=-4
\end{equation*}
    where $N_{b}:\mathcal{B}\to \mathcal{B}$ is a $s$-to-1 proper map and $m$ denotes the connectivity of $\mathcal{B}$. Solving this equation for the factors of $-4$ and the integer values of $s\geq 1$ and $m\geq 1$, the unique solution becomes $(s,m)=(3,1)$. In either case $\mathcal{B}$ is simply connected. Since there are no more free critical points in $\mathcal{Q}$, it follows that $\jul(N_{-e_3})$ is connected.
\end{proof}

\subsection{Configuration and absence of Herman rings}\label{sec:NoHR}

We begin by discussing the possible configuration of a Herman ring cycle under the hypothesis that $\LA$ is an arbitrary period lattice and $b\notin \{0,-e_1, -e_2, -e_3\}$. From Lemma~\ref{lem:propN} (2), the condition over $b$ implies that any pole of $N_b$ is a half-period lattice point, that is 
\begin{equation}\label{eq:polosmedios}
 \text{Poles}(N_b)=\frac 12 \LA\setminus\LA.
\end{equation}

If $H_i$ and $H_j$ are two distinct Herman rings in $\bbC$, we say that $H_i$ and $H_j$ exhibit a \textit{nested configuration} if
    \[H_i\subset B_j\qquad\text{or}\qquad H_j\subset B_i,\]
    where $B_i$ and $B_j$ denote the bounded component of $\Cc\setminus H_i$ and $\Cc\setminus H_j$, respectively.

\begin{proposition}[Nested configuration of Herman rings]\label{prop:onepole}
    Let $\LA$ be a given period lattice and consider
    \[N_b(z)=z-\frac{\WLA(z)+b}{\WLA'(z)},\qquad\text{with } b\notin\{0,-e_1, -e_2, -e_3\}.\]
    Assume $N_b$ has a $k$-cycle of Herman rings, ${\mathcal H}=\{H_1,\ldots, H_k\}$, for some $k\geq 1$. Then each bounded component of $\Cc\setminus H_j$
    %$B_j$
    contains exactly one and the same pole, and 
    %both $H_j$ and $B_j$ are 
    $H_j$ is symmetric with respect to that pole for all $j$. Furthermore, if $k\geq 2$, then the rings in the cycle exhibit a nested configuration among themselves.
\end{proposition}
\begin{proof}
    ~The maximum modulus principle implies that for at least one ring domain $H:=H_j$, its component $B:=B_j$ contains a pole $p$. Since $H$ surrounds $p$, then Proposition~\ref{prop:symm} implies that $H$ must be a symmetric Fatou component with respect to $p$, for otherwise, its symmetric image, $H'=2p-H$, would be a $k$-periodic Fatou component that surrounds $p$, and since $z\mapsto 2p-z$ is an isometry, then $H\cap H'\neq \emptyset$, implying their equality. Consequently, $B$ is also symmetric with respect to $p$.
    
    Assume $B$ contains two distinct poles, say $p_1:=p$ and $p_2$, then by symmetry with respect to $p_1$ the component $B$ would contain three distinct poles $p_1, p_2$ and $p_3=2p_1-p_2$. Now, reflecting $p_1$ and $p_3$ with respect to $p_2$ gives two more poles, $2p_2-p_1$ and $3p_2-2p_1$, which lie in $B$ along with $p_1,p_2,p_3$. Clearly, all five poles are distinct. Since every fundamental domain for the action of $\LA$ has exactly three poles, $B$ (and hence $H$) must lie outside a fundamental domain, which is impossible by Proposition~\ref{prop:bounded}. This implies that every $B_j$ contains at most one pole.
        
    Given $z,z'\in H$ with $z'=2p-z$, the condition $N_b^{m}(z')=2p-N_b^{m}(z)$ for all $m\geq 1$ implies that every Herman ring in $\mathcal H$ contains two points that are symmetric with respect to $p$, therefore, every ring is symmetric with respect to $p$. In particular, if $k\geq 2$ then each ring surrounds the pole $p$ and only this pole, so the $k$-cycle exhibits a nested configuration.
\end{proof}

One of the main steps in the proof of Proposition~\ref{prop:onepole} resides in the fact that the set of poles of the Newton's method is \textit{exactly} the set of half-periods $\frac 12 \LA\setminus\LA$, so each pole is a point of symmetry as described in Proposition~\ref{prop:symm}. It is easy to see that any other elliptic function of the form $\fLA=a(\WLA-b)^m$ with $a\in \bbC^*$, $b\notin\{e_1,e_2,e_3\}$ and $m\in\bbN$, also satisfies this hypothesis, so the nested configuration of Herman rings applies to $N_{\fLA}$.

\begin{remark}[Additional critical orbits and Herman rings.]\label{rem:cvs}
    If there exists a $k$-cycle of Herman rings, it is well known that $\partial \mathcal H:=\cup_{j=1}^k \partial H_j$ lies in the closure of the forward orbit of $\text{Sing}(N_b^{-1})$ (see, for example, \cite[Theorem 7]{Ber}). Thus, there exists at least one free (hence, additional) critical point, say $c$, whose forward orbit accumulates at $\partial \mathcal H$. Since $\mathcal H$ is symmetric with respect to a unique pole $p$, then $c'=2p-c$ is a second free critical point that also accumulates at $\partial \mathcal H$.
    
    If $\LA$ is triangular, then $\{c,c'\}+\LA$ represents all additional critical points of $N_b$. This, together with the $\LA$-equivariance property, imply that no forward orbit of a free critical point can enter any of the Herman rings, for otherwise there will be not enough free critical points to accumulate in $\partial \mathcal H$.
 \end{remark}   
The next result establishes the location of at least two free and symmetric critical values of $N_{b}$ (over any period lattice) with respect to the $k$-cycle of Herman rings.
\begin{lemma}\label{lema:cvalsin}
    Let $\LA$ be an arbitrary period lattice and assume $N_b=N_{\WLA,b}$ has a $k$-cycle of Herman rings for some $k\geq 1$. Denote by $H$ the innermost Herman ring in its $k$-cycle and let $B$ denote the bounded component in $\Cc \setminus H$. Then, $B$ contains at least two free critical values with disjoint residual classes. If $\LA$ is triangular, then $B$ contains exactly two free critical values and no other critical value.
\end{lemma}

\begin{proof}
    ~Let $\gamma$ be a curve in the foliation of $H$ disjoint from critical values and let $D=\text{int}(\gamma)$, which is bounded by definition. Propositions~\ref{prop:bounded} and \ref{prop:onepole} imply that $\overline{D}$ is contained in a fundamental parallelogram of $\LA$ and it is symmetric with respect to a unique pole $p$.
    First observe that $D$ cannot contain any superattracting fixed point: if it does, then by symmetry $D$ contains the closures of two superattracting immediate basins. On their boundaries lie two symmetric repelling fixed points, implying that $\gamma$ surrounds two distinct lattice points, a contradiction.
    
    Now assume $D$ is disjoint from $\sing(N_b^{-1})$. 
    Then $N_b^{-1}(\overline{D})$ consists of an infinite collection of topological disks, each of them mapping in a univalent fashion onto $\overline{D}$. Among these collection of preimages there must exist a topological disk $D'\subset N_b^{-1}(D)$, whose boundary is $N_b^{k-1}(\gamma)$. Then, $p\in D'$ and thus $N_b(D')$ is unbounded, contradicting the choice of $D$.

    Consequently, if $D$ contains one free critical value, say $v$, then it contains its symmetric image $v'=2p-v$, which is also a free critical value and it is disjoint from the residual class of $v$. If $\LA$ is triangular, then $D$ contains these two free critical values and no other by Remark~\ref{rem:cvs}.
\end{proof}

We are now ready to provide the proof of the first part of Theorem~\ref{thm:NoHR}. To achieve this, we restrict ourselves to the case of triangular lattices.

\begin{lemma}\label{lem:NoHR}
     Let $\LA$ be a triangular lattice. Then, for any $b\in \bbC$, the Newton map
    \[N_b(z) = z - \frac{\WLA(z)+b}{\WLA'(z)}\]
    does not have cycles of Herman rings of any period.
\end{lemma}
\begin{proof}
        ~If $b\in \{0,-e_1, -e_2, -e_3\}$ then Lemma~\ref{lem:bees} shows that $\jul(N_b)$ is connected, implying the absence of cycles of Herman rings.
        To obtain a contradiction, assume the existence of a parameter $b\notin\{0,-e_1, -e_2, -e_3\}$ for which $N_b$ has a $q$-cycle of Herman rings, $\mathcal{H}=\{H_1,\ldots,H_q\}$, with $q\geq 1$. Proposition~\ref{prop:onepole} implies that $\mathcal{H}$ is a nested cycle surrounding a unique pole $p\in \frac 12 \LA\setminus\LA$ 
        and each ring is symmetric with respect to $p$. Let $H_{in}$ and $H_{out}$ denote, respectively, the most inner and the most outer Herman rings in the nested configuration of the $q$-cycle (if $q=1$ the choice is obvious). Select any curve $\gamma_{in}\subset H_{in}$ in its foliation and denote by  $\Gamma$ the forward orbit of $\gamma_{in}$. 
        Let $\gamma_{out}$ be the unique curve in $\Gamma$ that lies in $H_{out}$.

        It follows from Proposition~\ref{prop:bounded} the existence of a fundamental parallelogram $\mathcal Q$, for which $\mathcal H$ (and hence $\Gamma$ and $p$) lies in the interior of $\mathcal Q$, while from the choice of the parameter, any other pole in $\mathcal{Q}$ lies in its boundary. 

        Since $\LA$ is triangular, let 
        $\mathcal{C}=\{c_1, c_2\}$ denote the set of free critical points in $\cQ$, where $c_2=2p-c_1$ and each $c_i$ is a double critical point. Set $v_i:=N_b(c_i)$ for each $i$. From Lemma~\ref{lema:cvalsin}, 
        we have that $v_1, v_2\in \text{int}(\gamma_{in})$.
        Define a simple curve $L_1\subset\bbC$ with extreme points at $v_1$ and $\infty$, avoiding all other critical values of $N_b$ (this is possible since $N_b(\Crit(N_b))$ is a discrete set). Let $L_2$ be its symmetric image through $p$, so that $L_2=2p-L_1$ as sets, and thus, also avoiding all other critical values except for $v_2$. From the assumption over $b$, all poles of $N_b$ are simple. Hence, the pullback of $L_1$ into $\mathcal Q$ under $N_b$ defines a graph, $T_1\subset \mathcal Q$, with  
        three 
        edges emanating from $c_1$ and each edge ending at a simple pole. By symmetry, $T_2$ is also a graph with the same properties. We show in Figure~\ref{fig:tripod} a sketch of the configuration of  $T_1\cup T_2$. 
\begin{figure}[htp]
    \centering
    \includegraphics[width=0.7\linewidth]{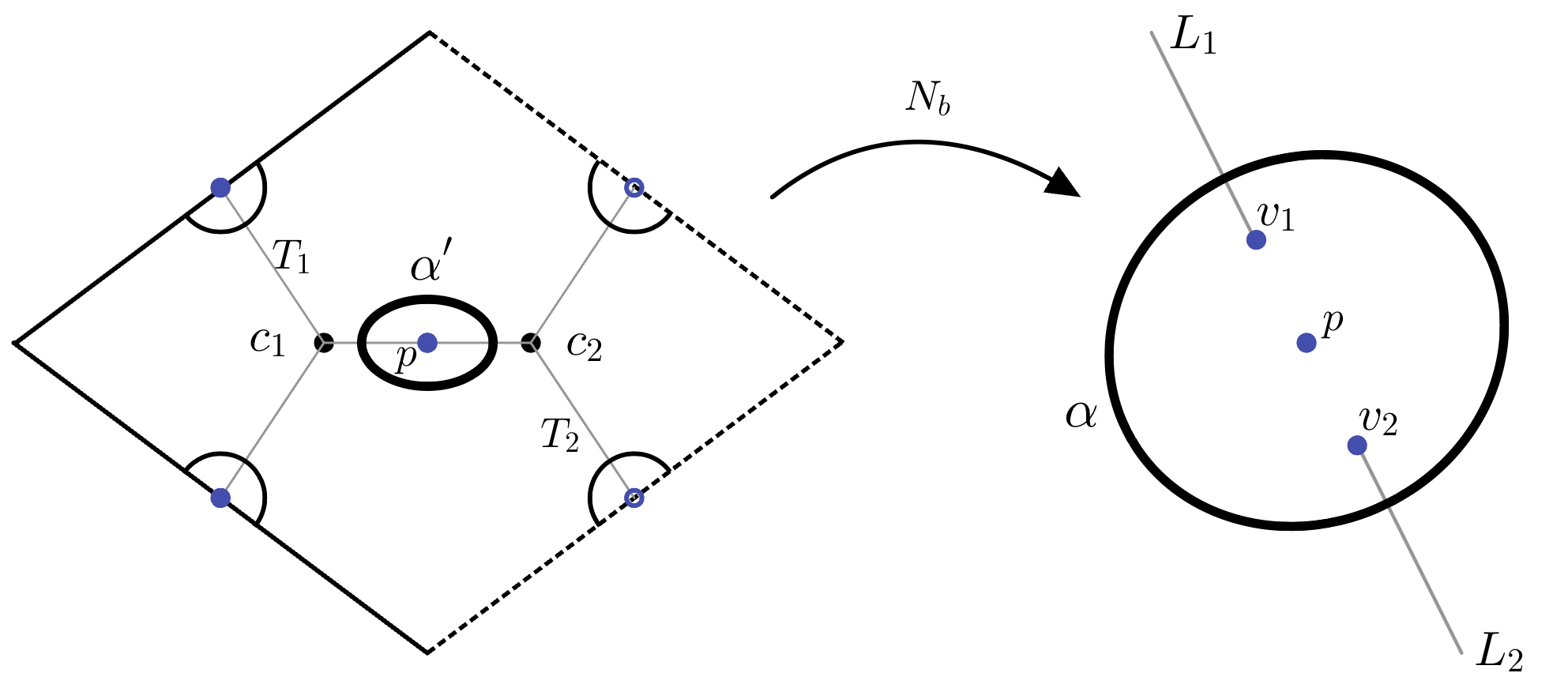}
    \caption{A sketch of the graphs $T_1$ and $T_2$ are shown in grey for the case $g_2=0$. Each $T_i$ connects a critical point $c_i$ with three poles of $N_b$, shown as blue dots, up to identification through a dashed side with its opposite side in $\mathcal{Q}$. The image of $T_1\cup T_2$ under $N_b$ is a curve $L_1\cup L_2$, symmetric with respect to $p$ and passing through $\infty$. Any curve $\alpha\in \Gamma$ has three preimages in $\mathcal{Q}$ up to identification, among them, $\alpha'$ is the only preimage surrounding $p$ and thus intersecting $T_1\cup T_2$ at the edge $[c_1, c_2]$.}
    \label{fig:tripod}
\end{figure}

    Since $v_1, v_2\in \text{int}(\gamma_{in})$, the curves $L_1$ and $L_2$ intersect every $\alpha\in \Gamma$. Deforming the curves $L_i$ if necessary, we can assume that the intersections of $L_1\cup L_2$ with each curve in $\Gamma$ are transversal and occur at two points, one in each $L_i$. For any $j=0,\ldots,q-1$, there exists an open neighbourhood $U_w$ of each $w\in N_b^j(\gamma_{in})\cap (L_1\cup L_2)$ where the map $N_b^{-1}|U_w$ is well defined and conformal. In conclusion, $T_1\cup T_2$ intersects each element in $\Gamma$ transversally and at two points. This can only be achieved if every $\alpha\in\Gamma$ intersects $T_1\cup T_2$ only at the edge $[c_1,c_2]$ containing $p$, yielding that $\mathcal{C}$ lies in the unbounded component $\Ext(\gamma_{out})$ and, in general, $\Crit(N_b)\subset \Ext(\gamma_{out})$.
    
We now proceed to show that $N_b$ cannot realise a $q$-cycle of Herman rings for any $q\geq 1$.
     
        If $q=1$, write $\gamma=\gamma_{in}=\gamma_{out}\subset H$. Then $\text{int}(\gamma)$ is disjoint from $\Crit(N_b)$ and contains a simple pole $p$. Thus, $\text{int}(\gamma)$ is mapped to $\Ext(\gamma)$ in a univalent manner under $N_b$. This would imply that any other foliation curve in $\text{int}(\gamma)\cap H$ is mapped to a different foliation curve in $\Ext(\gamma_{in})\cap H$, which is impossible since every foliation curve in $H$ is invariant under $N_b$.

For $q\geq 2$ let $A:=\text{int}(\gamma_{out})\cap \Ext(\gamma_{in})$ and observe that $\Crit(N_b)\subset \Ext(\gamma_{out})$. If $q= 2$, then $N_b|A:A\to A$ is univalent and bounded. Therefore, $\{N_b^m|A\}_{m\geq 1}$ is a normal family, which is impossible since $A\cap \jul(N_b)\neq \emptyset$.

        Finally, assume $q\geq 3$. Let $\alpha=N_b^{q-1}(\gamma_{in})\in \Gamma$ and observe that $\alpha$ is essentially contained in $A$. Then the annulus $B=\text{int}(\alpha)\cap \Ext(\gamma_{in})$ is also essentially contained in $A$. Moreover, $N_b|A$ is conformal and, by construction $N_b(B)=A$, thus $\text{mod}(B)=\text{mod}(A)$. Since $q\geq 3$ then $A$ admits a decomposition $A=B\sqcup (A\setminus \overline{B})$, where $A\setminus \overline{B}$ is an annular domain with non-empty intersection with at least $H_{out}$, implying a positive conformal modulus. At the same time, Gr\"otzsch's inequality leads to $\text{mod}(A\setminus \overline{B})=0$, a contradiction.        
\end{proof}

\subsection{Absence of unbounded Fatou domains}
As shown in Theorem~\ref{thm:unb_comp}, if an unbounded Fatou component is realised by the Newton's method of an odd or even elliptic function with arbitrary period lattice, such unbounded component must be simply connected and either a Baker domain or an escaping wandering domain. Firstly, we describe the topology of unbounded Fatou components.

%\textcolor{gray}{
Let $U_0\subset\fat(N_f)$ be an unbounded Fatou component of the Newton's method of an elliptic function $f$ with respect to a lattice $\LA$. Following \cite{HK24}, we say $U_0$ is a {\it toral band} if it contains an open subset $V$ which is simply connected in $\bbC$, but $V$ projects to a topological band around the torus $\TLA$ (see Remark \ref{rem:equivariance}) containing a homotopically nontrivial curve. It is called a {\it double toral band} if it projects to a set that contains closed paths that generate the fundamental group $\pi_1(\TLA)$. $U_0$ is called a {\it single toral band} if is not double, and {\it nonperiodic} if it maps onto a cycle of Fatou components but is not part of the cycle. It is not difficult to see that, from Theorem \ref{thm:unb_comp}, any toral band $U_0\subset\fat(N_f)$ is single.%} 

Theorem~\ref{thm:NoBaker} establishes the absence of unbounded Fatou components for the Newton's method of $\WLA+b$ over triangular lattices. Its proof requires the next lemma.
\begin{lemma}
    \label{lem:uni_unb}
    Fix a triangular lattice $\LA$ and let $N_b$ denote the Newton's method of $\WLA+b$. Suppose $U\subset\fat(N_b)$ is an unbounded component. Then $N_b$ is univalent in $U$ and its boundary is contained in the closure of the post-singular set. 
\end{lemma}
\begin{proof}
    ~From Theorem \ref{thm:unb_comp}, $U$ is simply connected. It follows that $U$ must be a {\it single toral band} for $N_b$. Let $W=\Pi(U)$, where $\Pi:\bbC\to\TLA$ is the canonical projection. Let $n_b
    :\TLA\setminus\Pi(\Poles(N_b))\to\TLA$ be the projection map as defined in equation (\ref{def:proj_map}). Then $W\subset \fat(n_b)$ is a 2-connected component. It follows from \cite{E93} that $n_b$ is a map of finite type, and from the invariance of the Fatou set of $n_b$, we may assume without loss of generality that $W$ is invariant under $n_b$. 
    From the Riemann-Hurwitz formula, we have that $n_b$ is univalent on $W$, which implies $N_b$ is univalent on $U$. Moreover, $W$ is a Herman ring for $n_b$, then $\partial W$ is contained in the post-singular set.  
\end{proof}

\begin{lemma}\label{thm:NoBaker}
Fix a triangular lattice $\LA$ and let $N_b$ denote the Newton's method of $\WLA+b$. Then, for any $b\in\bbC$, $N_b$ has no unbounded Fatou components. In particular, $N_b$ has no Baker domains.
\end{lemma}

\begin{proof}
    ~Suppose $U\subset \fat(N_b)$ is an unbounded component. We know $U$ is simply connected and $N_b$ is univalent on $U$. From Lemma \ref{lem::inv-LA}, we may assume $U$ intersects the principal fundamental parallelogram $\cQ_0$. 
    From Lemma \ref{lem:uni_unb}, the free critical point $c_1\in\cQ_0$ is such that $c_1\in\jul(N_b)$. Let $c_2=2\omega_3-c_1\in\cQ_0$ 
    the other critical point, symmetric with respect to $\omega_3\in\cQ_0$. 
    Let $v_i=N_b(c_i)\in\jul(N_b)$ be the corresponding critical values and, as in the proof of Lemma~\ref{lem:NoHR}, consider the symmetrical curves $L_i$ connecting $v_i$ with the point at infinity. For this proof, assume also that $L_i\cap V=\emptyset$ for any unbounded Fatou component $V$. If $T\subset\overline{\cQ_0}$ denotes the union of the tripods $T_1$ and $T_2$, then $T$ connects
    opposite sides of the parallelogram $\cQ_0$, see Figure \ref{fig:tripod}. 
    
    Hence, $U$ must intersects one of the 
    sets $T_i$, let $I=T\cap U$. Then, $N_b(I)$ (and hence $N_b(U)$) is contained in a bounded component. This is a contradiction, since $N_b$ has no finite asymptotic values.
\end{proof}
The proof of Theorem~\ref{thm:NoHR} is now a consequence of   Lemma~\ref{lem:NoHR}, Lemma~\ref{thm:NoBaker}, and Theorem \ref{thm:conn_gen}.

\section{Existence of wandering domains}\label{sec:wan}
Here we present the proof of Theorem \ref{thm:wandering}, concerning the existence of wandering domains for the Newton's method of an odd or an even elliptic function with an arbitrary period lattice. For triangular period lattice $\LA$, we explicitly describe in Proposition~\ref{ex:wanderN_b} the parameters $b\in \bbC$ for which the Newton's method of $\WLA+b$ has bounded and simply connected wandering domains. In both results, wandering domains coexist with invariant attracting basins in the same dynamical plane.

The proof of Theorem~\ref{thm:wandering} relies on an auxiliar transcendental function: if $\la\in \LA^*$ is a given non-zero period, let
\[F_\la(z):=N_{\fLA}(z)+\la,\]
where $N_{\fLA}$ is the Newton's method for a given elliptic function $\fLA\in\ALA$. Observe that $F_{\la}$ and $N_{\fLA}$ are transcendental meromorphic functions that commute: Indeed, if $T_\la(z)=z+\la$, then
    \[F_{\la}\circ N_{\fLA}=T_\la\circ N_{\fLA}\circ N_{\fLA}=N_{\fLA}\circ T_\la \circ N_{\fLA}=N_{\fLA}\circ F_{\la},\]
    whenever both sides are defined. Since $F_{\la}$ is a constant translation of $N_{\fLA}$, they share the same sets of critical points and poles. As the set of prepoles has infinite cardinality,  
    Theorem 1.5 in \cite{Tsa} implies that $\jul(F_{\la})=\jul(N_{\fLA})$, therefore, their Fatou sets coincide. 
\vspace{0.5cm}

\noindent{\bf Theorem D.}
For any given lattice $\LA$, let $N_{\fLA}$ denote the Newton's method for an odd or even elliptic function $\fLA$. Consider the function $F_{\la}(z)=N_{\fLA}(z)+\la$ for some fixed $\la\in \LA^*$ and let $q\geq 1$. If $F_{\la}$ has a $q$-cycle of (super)attracting  basins containing only free critical points,  %in $\Zeros(\fLA'')\setminus \Zeros(\fLA)$},
then $N_{\fLA}$ has (at least) $2q$ distinct orbits of escaping wandering domains that coexist with the invariant (super)attracting basins for the roots of $\fLA$.

\begin{proof}
    ~Fix $\la\in\LA^*$ and assume that $F_{\la}$ has a (super)attracting $q$-cycle   
    $C=\{z_0, z_1,\ldots,z_{q-1}\}$, 
labeled in such a way that $c_0$, a free critical point in $\Zeros(\fLA'')$, %\setminus \Zeros(\fLA)$}\footnote{\textcolor{blue}{no sería suficiente decir ``free critical point''?}},
lies in the immediate basin  of $z_0$. From the $\LA$-equivariance of $N_{\fLA}$,
    then $F_{\la}^2(z)= N_{\fLA}(N_{\fLA}(z)+\la)+\la=N_{\fLA}^2(z)+2\la$, and in general $F_{\la}^k(z)=N_{\fLA}^k(z)+k\la$ for all $k\in \bbN$ and all $z$ away from the prepoles of $F_{\la}$, which coincide with the
    prepoles of $N_{\fLA}$. Using the periodicity of $z_0$ under $F_{\la}$, the orbit of $z_0$ under $N_{\fLA}$ becomes
    \begin{equation}\label{eq:iterN}
    N_{\fLA}^{m}(z_0)=
    F_{\la}^j(z_0)-m\la
    \end{equation}
    for all $m\geq 0$ and  $j=j(m)\in\{0,\ldots,q-1\}$ so that $q|(m-j)$.

    Denote by $U_0,\ldots,U_{q-1}$ the pairwise disjoint Fatou components for $F_{\la}$ that form the immediate basin of $C$, with $z_j\in U_j$ for $j=0,\ldots, q-1$. Then $U_0,\ldots, U_{q-1}$ are also pairwise disjoint Fatou components for $N_{\fLA}$, and the identity in (\ref{eq:iterN}) implies that 
    \[N_{\fLA}^{m}(U_0)=U_j-m\la.\]
    
    We claim that any positive iterate of $U_0$ under $N_{\fLA}$ belongs to the immediate basin of a (super)attracting $q$-cycle under $F_{\la}$, distinct from $C$. First, observe that $F_{\la}$ commutes with $T_{\la'}$ for any $\la'\in \LA^*$, as
    \[T_{\la'}\circ F_{\la}=T_{\la'}\circ (N_{\fLA}\circ T_\la)=N_{\fLA}\circ T_{\la'}\circ T_{\la}=N_{\fLA}\circ T_\la\circ T_{\la'}=F_{\la}\circ T_{\la'}.\]
    Therefore, any critical point of the form $c_0+\la'$ lies in a (super)attracting $q$-cycle under $F_{\la}$ given by $C+\la'=\{z_0+\la', z_1+\la',\ldots, z_{q-1}+\la'\}$, with immediate basin formed by $U_0+\la', U_1+\la',\ldots, U_{q-1}+\la'$, and $c_0+\la'\in U_0+\la'$. 
    
    Observe that $C$ and $C+\la'$ are disjoint, for if $z_i=z_0+\la'$, and $j>0$ is such that $i+j=n$, then $z_{i+1}=F_\la(z_i)=z_1+\la'$, and (inductively) $z_{i+j}=z_q=z_0=z_j+\la'$, which implies that $z_i=z_0+\la'=(z_j+\la')+\la'=z_j+2\la'$. Now, $z_0=z_{i+j}=z_j+\la'$, obtaining $z_i=z_{2j}+3\la'$. Inductively, $z_i=z_{(k-1)j}+k\la'$, which contradicts the cyclic nature. Therefore, $C$ and $C+\la'$ are disjoint.
    Therefore $U_0\neq U_j+\la'$ for all $\la'\in \LA^*$ and all $j\in\{0,\ldots,q-1\}$. 
    
   To show that $U_0$ is a wandering Fatou component for $N_{\fLA}$, assume there are distinct integers $m, \ell\geq 0$ so that $N_{\fLA}^m(U_0)=N_{\fLA}^\ell(U_0)$. That is, there exist $j,i\in\{0,\ldots,q-1\}$ such that $U_j-m\la=U_i-\ell\la$. If $i\geq j$, then, after applying a translation followed by  $F_{\la}^{q-j}$, this is equivalent to
   \[U_0=U_{i-j}-(\ell-m)\la\]
    which is impossible by the previous argument. The case $i\leq j$ is similar. By the same kind of reasoning it follows that the orbit of $U_j$ under $N_{\fLA}$ is wandering for each $j\in\{0,\ldots,q-1\}$, obtaining $q$ distinct orbits of %{\bf \sout{simply connected}} 
    wandering domains. The $\LA$-equivariance implies that these wandering domains are necessarily of escaping type. At the same time, the forward orbits of $U_j$ under $N_{\fLA}$ are clearly disjoint from the invariant basins of the roots $\Zeros(\fLA)$, therefore $\fat(N_{\fLA})$ contains invariant attracting basins and wandering domains.

    There is a second collection of $q$ distinct orbits under $N_{\fLA}$ of %{\bf \sout{simply connected}} 
    wandering domains that arise from applying the above arguments to $F_{-\la}$. Indeed, let $b$ and $\la$ be as before. Since $c_1=\la-c_0$, then
    \[F^k_{-\la}(c_1)=N_{\fLA}^k(\la-c_0)-k\la=\la-(N_{\fLA}^k(c_0)+k\la)=\la-F_{\la}^k(c_0),\]
    so the orbit of $c_1$ under $F_{-\la}$ is either $q$-periodic or converges to an attracting $q$-cycle, $C'=\{w_0,\ldots,w_{q-1}\}$, with $w_j=\la-z_j$. That is, $F_{-\la}$ has a (super)attracting $q$-cycle with a free critical point. The rest of the arguments are essentially the same but applied to $F_{-\la}$.    
\end{proof}

\begin{remark}
    If $U\subset\fat(N_{\fLA})$ is one of the $2q$ escaping wandering domains given by Theorem~\ref{thm:wandering}, the $\LA$-equivariance property of $N_{\fLA}$ implies that $U+\la'$ is also a wandering domain for any $\la'\in \LA$, therefore $\fat(N_{\fLA})$ consists of infinitely many wandering domains. Among them, there are $2q$ distinct \emph{orbit behaviours}: if $U$ escapes under iteration along a given pattern (say, $N_{\fLA}(U)\to \infty$ along the real axis) then $U+\la'$ also escapes along the same pattern (say, parallel to the real axis).
\end{remark}

\begin{figure}[htp]
    \centering
    \includegraphics[width=0.9\linewidth]{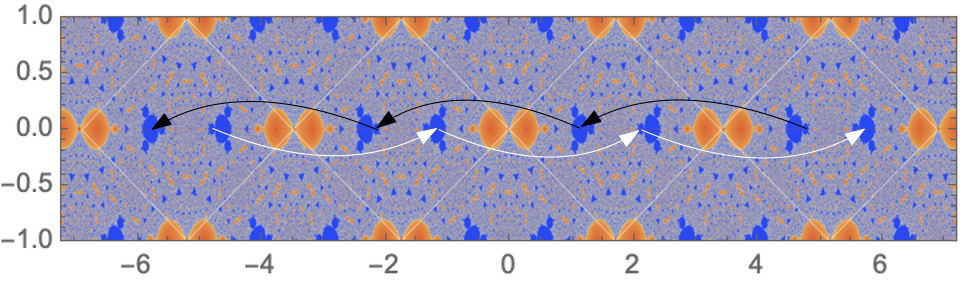}
    \caption{The dynamical plane for the Newton method of $\WLA(z)+b$, where $\LA$ is triangular (in horizontal position) with invariants $g_2=0$ and $g_3\approx -12.8254$. As explained in Proposition~\ref{ex:wanderN_b}, the parameter $b\approx -11.68$. The superattracting basins of the roots of $\WLA(z)+b$ are shown in orange while two cycles of simply connected wandering domains are displayed in blue, with black and white arrows indicating their escaping orbits. The $\LA$-equivariance property  implies the existence of infinitely many wandering domains with essentially two distinct orbital behaviours; one escaping to the left and one escaping to the right.} %of each cycle.}

    \label{fig:wandering-arrows}
\end{figure}

\begin{proposition}\label{ex:wanderN_b}       Let $\LA$ be a triangular lattice. For any $\la\in \LA^*$ and for any $c\in \WLA^{-1}(0)$, if $b=\la \WLA'(c)$ then the Fatou set $\fat(N_b)$ consists of simply connected invariant attracting basins %for the zeros of $\WLA+b$, and %two classes of 
      and simply connected escaping wandering domains. %{\bf \sout{Moreover, $\jul(N_b)$ is connected.}}
  
%      Let $\LA$ be a triangular lattice. For any $\la\in \LA^*$ and for any $c\in \WLA^{-1}(0)$, if $b=\la \WLA'(c)$ then the Newton map $N_b(z)$ has two distinct orbits of simply connected escaping wandering domains. Moreover, $\jul(N_b)$ is connected.
\end{proposition}

\begin{proof}~Fix $\la\in \LA^*$, select $c\in \WLA^{-1}(0)$ and set $b=\la \WLA'(c)$. Since $\LA$ is triangular, Theorem~\ref{thm:NoHR} shows $\jul(N_b)$ is connected, so any Fatou component is simply connected. Moreover, the shape of the lattice also implies that $\WLA'(c)\neq 0$, therefore $b\neq 0$ and the roots of $\WLA+b$ generate superattracting invariant basins for $N_b$; observe that $c$, a free critical point, may or may not be contained in one of these basins. Nonetheless, the auxiliar map $F_{\la}(z)=N_b(z)+\la$ has a superattracting fixed point at $c$, as the choice of $b$ implies that
\[N_b(c)=c-\frac{0+\la\WLA'(c)}{\WLA'(c)}=c-\la,\qquad\text{thus}\qquad F_{\la}(c)=N_b(c)+\la=c,\]
and $F_\la'(c)=N_b'(c)=0$, therefore $c$ is a superattracting fixed point of $F_\la$. Denote by $U_c\subset\fat(F_\la)$ the immediate basin of $c$.

%, which is simply connected by Theorem~\ref{thm:NoHR}.
%since no other (free) critical points lie in $U_c$. Indeed, as shown in the proof of Theorem~\ref{thm:wandering}, $F_{\la}$ has  superattracting fixed points at \textcolor{blue}{$\{c\}+\LA$} while $F_{-\la}$ has superattracting fixed points at \textcolor{blue}{$\{\la-c\}+\LA$}, and the functions $N_b, F_{\la}$ and $F_{-\la}$ commute pairwise. \textcolor{blue}{Since $\LA$ is triangular, there are only two residual classes of free critical points given by $[c]$ and $[\la-c]$, and therefore}, each free critical point of $N_b$ is contained in a superattracting invariant basin \textcolor{blue}{of $F_\la$ or $F_{-\la}$} and no other free critical point lies in that component. This, combined with the fact that each superattracting fixed point in $\WLA^{-1}(-b)$ lies in its own invariant component, implies that every Fatou component is simply connected and consequently $\jul(N_b)$ is connected.\footnote{\textcolor{blue}{la conectividad se puede probar directamente con el Teorema A, o no?}}

The maps $N_b$ and $T_\la$ commute and $N_b(c)=c-\la$, therefore $N_b^k(c)=c-k\la$ for all $k\geq 1$.  
The images of $U_c$ under $N_b$ are given by $N_b^k(U_c)=U_c-k\la$, thus $N_b$ acts as a translation over $U_c$ and consequently $U_c$ is a simply connected escaping wandering domain for $N_b$. The symmetric image $U'=-U_c+\la$ is invariant for $F_{-\la}(z)=N_b(z)-\la$ and thus it is also wandering under $N_b$ since
\[N_b^k(U')=N_b^k(\la-U_c)=-N_b^k(U_c)+\la=-U_c+(k+1)\la=U'+k\la.\]
\end{proof}

%%
%\textcolor{blue}{
%\begin{proposition}%\label{ex:wanderN_b}   
%      Let $\LA$ be a triangular lattice. For any $\la\in \LA^*$ and for any $c\in \WLA^{-1}(0)$, if $b=\la \WLA'(c)$ then the Fatou set $\fat(N_b)$ consists of invariant attracting basins %for the zeros of $\WLA+b$, and %two classes of       and simply connected escaping wandering domains. Moreover, $\jul(N_b)$ is connected.
%\end{proposition}
%\begin{proof}
%    (\textbf{Sketch of additional argument:} $\LA$ triangular implies:
%    \begin{itemize}
 %       \item The zeros are super-attracting fixed points of $N_b$. 
 %       \item $\LA$ the only additional critical points are the zeros of $\WLA$, that must be associated to the additional super-attracting point of $F_\lambda$.
 %       \item Every fundamental parallelogram contains a member of the classes obtained in Theorem D. 
%    \end{itemize}
%\end{proof}
%}\footnote{una alternativa}

\begin{figure}[htp]
\centering
\begin{subfigure}[b]{0.47\linewidth}
\includegraphics[width=\linewidth]{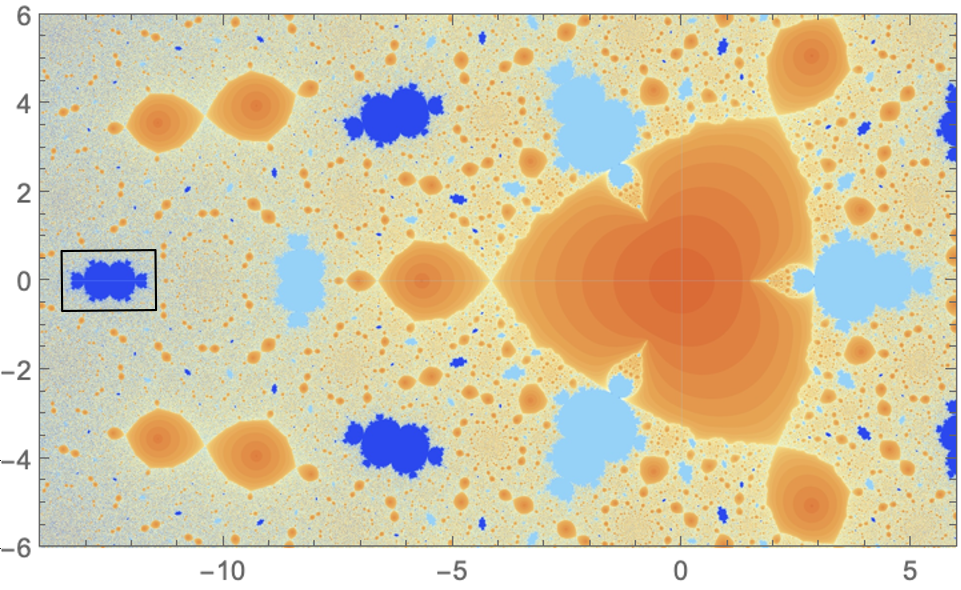}
\caption{}
\label{fig:Tpara_space}
\end{subfigure}
\begin{subfigure}[b]{0.45\linewidth}
\includegraphics[width=0.9\linewidth]{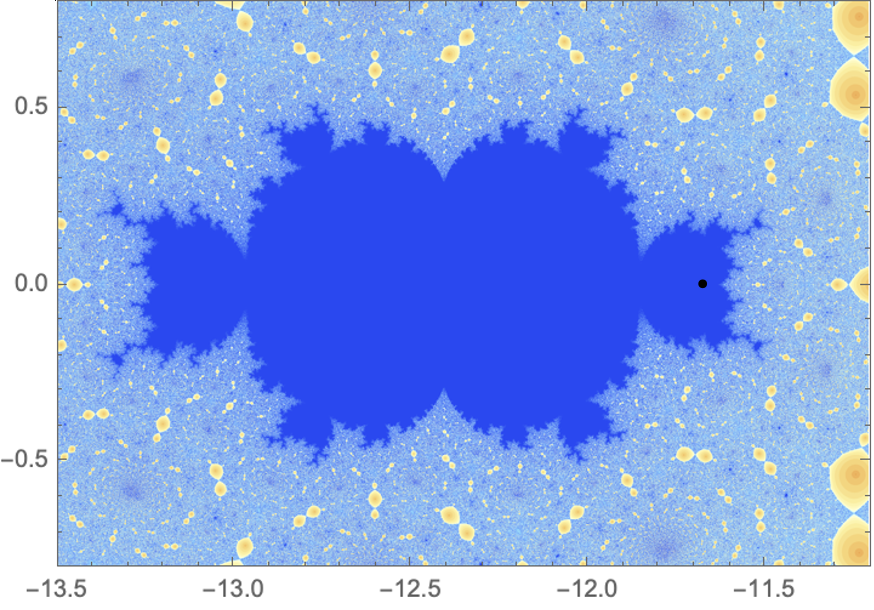}
\caption{}
\label{fig:Zpara_space}
\end{subfigure}
\caption{(a) The $b$-parameter space for the Newton's method of $\WLA(z)+b$ with $\LA$ triangular. Orange regions consist of parameters for which the free critical points of $N_b$ are trapped by the attracting basins of roots. Blue multibrots correspond to regions of parameters that give rise to wandering domains while light-blue multibrots represent parameters for which a free critical points belong to the basin of 
periodic cycles. (b) A magnified view of a ``wandering'' Multibrot inside the black rectangle in (a) is shown, with parameter $b\approx -11.68$ marked with a black dot.
}
\label{fig:para_space}
\end{figure}

As was mentioned before, Theorem~\ref{thm:wandering} and, in particular, Proposition~\ref{ex:wanderN_b}, provide the first examples of wandering domains for Newton's method of elliptic functions coexisting with the attracting basins of the zeros of the function. Moreover, it is not difficult to see that, from the proof of Proposition~\ref{ex:wanderN_b}, wandering domains for $N_b$ are preserved under any small perturbation of the parameter $b$. In other words, there is an open set of parameters for which the Newton's method $N_b$ has wandering domains, see Figure \ref{fig:para_space} and compare with the results in \cite{FF25}. We remark that the numerical experiments shown in Figure~\ref{fig:wandering-arrows} have been performed for a triangular lattice in horizontal position with $\lau=\exp(\pi i/6)$, $\lad=\bar{\lambda}_1$ and thus $g_3\approx -12.8254$.

A detailed study of the parameter space of the Newton's method
\[N_{\WLA,b}(z)=z-\frac{\WLA(z)+b}{\WLA'(z)},\]
with $\LA$ triangular, will be presented in an upcoming  article.

\bibliographystyle{amsplain}

\end{document}